\documentclass[12pt, oneside]{amsart}
\usepackage[utf8]{inputenc}

\usepackage[dvips]{geometry}
\usepackage[table]{xcolor}

\usepackage{arydshln}
\usepackage{graphicx}
\usepackage{amssymb}
\usepackage[section]{placeins}
\usepackage{todonotes}%%% this is to write comments%%%%%%

\usepackage{enumitem,amssymb}
\newlist{todolist}{itemize}{2}
\setlist[todolist]{label=$\square$}

\usepackage{amsfonts}
\usepackage{amsmath}
\usepackage{amsthm}
\usepackage{euscript}
\usepackage{amssymb}
\usepackage{mathrsfs}
\usepackage{tikz}
\usepackage{caption}
\usepackage{subfig}
\usepackage{mathtools}
\usepackage{nicematrix}
\usetikzlibrary{fit}
\usepackage[colorlinks]{hyperref}

\theoremstyle{plain}
\newtheorem{thm}{Theorem}[section]
\newtheorem{cor}[thm]{Corollary}
\newtheorem{lem}[thm]{Lemma}
\newtheorem{prop}[thm]{Proposition}

\theoremstyle{definition}
\newtheorem{defn}[thm]{Definition}
\theoremstyle{remark}
\newtheorem{rem}[thm]{Remark}
\newtheorem{ex}[thm]{Example}
\numberwithin{equation}{section}

\newcommand{\OP}{\mathcal{O}_{\mathbb{P}^1}}

\newcommand{\C}{\mathbb C}

\newcommand{\GG}{\mathbb G}
\newcommand{\PP}{{\mathbb P}}

\newcommand{\OO}{\mathcal{O}}

\newcommand{\LL}{\mathcal{L}}

\renewcommand{\ge}{\geqslant}
\renewcommand{\le}{\leqslant}

\DeclareMathOperator{\hh}{h}

\DeclareMathOperator{\rk}{rk}

\DeclareMathOperator{\GL}{GL}

\begin{document}

\title{Pencils of singular quadrics of constant rank and their orbits}

\author{Ada Boralevi}
\address{Dipartimento di Scienze Matematiche \lq\lq G. L. Lagrange\rq\rq, Politecnico di Torino, Corso Duca degli Abruzzi 24, 10129 Torino, Italy}
\email{\href{mailto:ada.boralevi@polito.it}{ada.boralevi@polito.it}}

\author{Emilia Mezzetti}
\address{Dipartimento di Matematica e Geoscienze, Sezione di Matematica e Informatica,  Universit\`a degli Studi di Trieste, Via Valerio 12/1, 34127 Trieste, Italy}
\email{\href{mailto:mezzette@units.it}{mezzette@units.it}}
\thanks{The first named author has been partially supported by MIUR grant Dipartimenti di Eccellenza 2018-2022 (E11G18000350001). The second named author has been  partially supported by the grant PRIN $2017$SSNZAW$\_005$ \lq\lq Moduli Theory and Birational Classification''. Both authors are members of INdAM--GNSAGA}

\subjclass[2020]{14C21, 14M15, 14L30, 14F06}

\keywords{Symmetric matrices, pencils of quadrics, general linear group, orbit}

\begin{abstract} 
We give a geometric description of singular pencils of quadrics of constant rank, relating them to the splitting type of some naturally associated vector bundles on $\PP^1$. Then we study their orbits in the Grassmannian of lines, under the natural action of the general linear group.
\end{abstract}

%\date{\today}

\dedicatory{Dedicated to Giorgio Ottaviani on the occasion of his $60^{th}$ birthday}

\maketitle
\sloppy

%%%%%%%%%%%%%%%%%%%%%%%%%%%%%%%%%%%%%%%%%%%%%%%%%%%%%%%%%%%%%%%%%%%%%%%%
%%%%%%%%%%%%%%%%%%%%%%%%%%%%%%%%%%%%%%%%%%%%%%%%%%%%%%%%%%%%%%%%%%%%%%%%
\section{Introduction}
%%%%%%%%%%%%%%%%%%%%%%%%%%%%%%%%%%%%%%%%%%%%%%%%%%%%%%%%%%%%%%%%%%%%%%%%
%%%%%%%%%%%%%%%%%%%%%%%%%%%%%%%%%%%%%%%%%%%%%%%%%%%%%%%%%%%%%%%%%%%%%%%%

A {\em pencil of quadrics} in the projective space of dimension $N$ is a two-dimensional linear subspace $\LL$ in the space of symmetric matrices of order $N+1$, and it is a widely studied object in algebraic geometry.

A complete classification of pencils of quadrics, based on algebraic considerations, Segre symbols and minimal indices, has been known for a long time: we refer to the classical book by Gantmacher \cite{Gantmacher} and the expository article by Thompson \cite{Thompson}. 

There is also an extensive literature on geometric descriptions and interpretations of pencils of quadrics; among the many contributions, let us cite some older works, from \cite{Segre} to \cite{Dimca}, as well as more recent ones, such as \cite{Fevola-Sturmfels}.

Often, when studying pencils of quadrics in $\PP^N$, one assumes that they are {\em regular}, that is, that they contain quadrics of maximal rank $N + 1$. As observed in \cite{Fevola-Sturmfels}, these pencils form an open subset in the appropriate Grassmannian, that admits a natural stratification by Segre symbols. The pencils in the complementary closed subset, called {\em singular pencils}, are less studied, even if in \cite{Gantmacher} it is shown that their analysis can be traced back to that of regular pencils and of singular pencils of constant rank. The purpose of this article is to give a description of the geometry of such pencils of constant rank, to relate it to the splitting of certain bundles on $\PP^1$ naturally associated with them, and to give a description of their orbits under the natural action of the general linear group $\GL(N+1)$.

To be more precise, we set up our notations: we work over an algebraically closed field of characteristic $0$, for simplicity over the complex field $\C$.
Let $V$ be a vector space of dimension $N+1$ over $\C$. Denote by $X$ the Veronese variety, that is, the image of the Veronese map $\PP(V) \to \PP(S^2V)$. The natural action of the group $\GL(N + 1)$ on $\PP(V)$ extends to $\PP(S^2V)$, and the orbits under this latter action are $X$ and its secant varieties. 

Fixing a basis for $V$, the elements of $S^2V$ can be seen as symmetric $(N+1) \times (N+1)$ matrices: then the action of $\GL(N+1)$ is the congruence, $X$ corresponds to symmetric matrices of rank 1, and its $k$-secant variety $\sigma_k(X)$ to symmetric matrices of rank at most $k$. 

\smallskip

Working in this projective setting, we interpret a pencil of quadrics as a line $\PP(\LL) \subseteq \PP(S^2V)$: it is {\em singular} when it is entirely contained in the determinantal hypersurface $\sigma_N(X)$. If a singular pencil is entirely contained in a stratum 
$\sigma_k(X) \setminus \sigma_{k-1}(X)$, we say that the pencil has {\em constant rank} $k$.
All the quadrics in such a pencil are cones having as vertex a linear space of dimension $N-k$.
\smallskip

In Section \ref{first results} we show that a pencil of constant rank $k$ corresponds to a matrix of linear forms in two variables, that naturally defines a map of vector bundles of rank $N+1$ over $\PP^1$; since the rank is constant, the  cokernel $E$ of this map is also a vector bundle over $\PP^1$, of rank $N+1-k$, and its first Chern class is $\frac{k}{2}$; in particular the constant rank $k$ is an even number that we denote by $2r$. We prove that the splitting type of $E$ characterizes the orbits, and for each orbit we give two explicit constructions for the canonical form of the representative: one is the expression described in \cite{Gantmacher}, the other one  is  analogous to the representative given in \cite{Fania_Mezzetti}, adapted from the skew-symmetric case. Indeed, several techniques used in articles on spaces of skew-symmetric matrices of constant rank, such as \cite{Manivel_Mezzetti, bo_me_piani,bor-fan-mez-quadriche}, can be applied to pencils of quadrics.

\smallskip

Analyzing these canonical forms, in Section \ref{examples} we describe the geometry of the pencils in the various orbits. If we make the assumption that the bundle $E$ has no trivial direct summand, which is equivalent to the condition that the quadrics in the pencil $\LL$ have no common point in their vertices, the pencil is called non-degenerate. In this case, if the splitting type of $E$ is $r_1, \ldots, r_h$, any two quadrics of $\LL$ have a generating space $S$ of (maximal) dimension $N-r$ in common, and are tangent along a rational normal scroll of dimension $r$ and type $r_1, \ldots, r_h$ contained in $S$.

In Section \ref{dimensione orbite}, we prove our main result Theorem \ref{main}: we find an explicit expression for the dimension of every $\GL(N+1)$-orbit of pencils of constant rank. We recall that these pencils are all unstable, nevertheless we are able to find an explicit expression for the matrices in the Lie algebra of the stabilizer of any pencil $\LL$. In particular these Lie algebras all have dimension $5$ when the corank of the pencil is $1$, i.e. $E$ is a line bundle with $c_1=r$. In Proposition \ref{descrizione algebra} we prove that they are of the form  $\mathfrak{sl}_2 \ltimes \C^2$. We conclude with a table collecting the results for $r\le 6$.

%%%%%%%%%%%%%%%%%%%%%%%%%%%%%%%%%%%%%%%%%%%%%%%%%%%%%%%%%%%%%%%%%%%%%%%%
%%%%%%%%%%%%%%%%%%%%%%%%%%%%%%%%%%%%%%%%%%%%%%%%%%%%%%%%%%%%%%%%%%%%%%%%
\section{Classification's details and first results}\label{first results}
%%%%%%%%%%%%%%%%%%%%%%%%%%%%%%%%%%%%%%%%%%%%%%%%%%%%%%%%%%%%%%%%%%%%%%%%
%%%%%%%%%%%%%%%%%%%%%%%%%%%%%%%%%%%%%%%%%%%%%%%%%%%%%%%%%%%%%%%%%%%%%%%%

Recall from the Introduction that, given an $(N+1)$-dimensional vector space $V$, one has the natural Veronese map $\PP(V) \to \PP(S^2V)$ sending $[v] \mapsto [v^2]$, whose image is the Veronese variety $X$. Once we fix a basis of $V$, the elements of $S^2V$ are identified with symmetric $(N+1) \times (N+1)$ matrices, $X$ corresponds to symmetric matrices of rank 1, and its $k$-secant variety $\sigma_k(X)$ to symmetric matrices of rank at most $k$. The group $\GL(N + 1)$ acts by congruence on 
$\PP(S^2V)$, and the orbits are exactly $X$ and its secant varieties. 

\smallskip

Now let  $\PP(\LL) \subseteq \sigma_k(X) \setminus \sigma_{k-1}(X)$ be a singular pencil of quadrics of constant rank $k$. Notice that $\PP(\LL)$ can be seen as a symmetric matrix whose entries are linear forms in two variables, that is, a vector bundle map on 
$\PP^1=\PP(\LL)$ of the form $V^* \otimes \OO_{\PP^1}(-1)\rightarrow V \otimes \OO_{\PP^1}$, inducing a long exact sequence:
\begin{equation}\label{succ esatta lunga}
0 \to E^*(-1) \to V^* \otimes \OO_{\PP^1}(-1) \to V \otimes \OO_{\PP^1} \to E \to 0.
\end{equation}
 
The cokernel  is a vector bundle of rank $N+1-k$ on $\PP^1$, hence it splits as a direct sum of line bundles; we denote it by $E$. The symmetry implies that the kernel is $E^*(-1)$.

From a direct computation of invariants (see \cite{Ilic_JM} for details), one finds that the rank $k=2r$ is even, the bundle $E$ is generated by its global sections, and moreover its first Chern class is $c_1(E)=r$. 

\smallskip

We start our description of $\LL$ generalizing to the symmetric case some results from \cite{Fania_Mezzetti} that refer to the skew-symmetric case. We are of course interested in non-trivial cases: for this, recall that a space of matrices is called {\em non-degenerate} if the kernels of its elements intersect in the zero subspace and the images of its elements generate the entire vector space $V$. This is equivalent to saying that the space is not $\GL(N + 1)$-equivalent to a space of matrices with a row or a column of zeroes. Therefore the classification of degenerate spaces of matrices can be traced back to that of non-degenerate spaces of matrices of smaller size. From now on, we will only consider non-degenerate spaces of constant rank $2r$.

Non-degeneracy also implies that as the quadrics vary in the pencil $\LL$, their vertices are pairwise disjoint. 

\smallskip

An immediate remark is that not all values of $N$ allow a non-degenerate pencil of symmetric matrices of size $N+1$ and fixed constant rank $2r$.

\begin{prop}
Let  $\LL \subset \PP(S^2V)$ be a non-degenerate pencil of singular quadrics of constant rank $2r$. Then $2r \le N \le 3r-1$.
\end{prop}

\begin{proof}
The proof of \cite[Proposition 3.6]{Fania_Mezzetti} goes through step by step. Since the cokernel bundle $E$ from \eqref{succ esatta lunga} is a vector bundle on $\PP^1$, it is of the form
\[E=\OP^{m_0} \oplus \OP(1)^{m_1} \oplus \cdots \oplus \OP(k)^{m_k},\]
where $m_0,\ldots,m_k$ are non-negative integers such that $m_1+2m_2 + \ldots + km_k=c_1(E)=r$, and $m_0+m_1+\ldots+m_k=\rk(E)=N+1-2r$.

The assumption that $\LL$ is non-degenerate implies $m_0 = 0$.

Obviously $2r \le N$. For the other inequality, we compute the cohomology of the sequence \eqref{succ esatta lunga} and deduce that $N+1 = \hh^0(V \otimes  \OP) \ge \hh^0(E)$, and 
\[\hh^0(E)= 2 \rk(E) + c_1(E)= 2(N+1-2r)+r. \qedhere\]
\end{proof}

\smallskip

The group $GL(N+1)$ acts by congruence on $\mathbb P(S^2V)$, the space of quadrics in $\PP(V)$, and thus it acts on pencils of quadrics, that correspond to lines in $\mathbb P(S^2V)$: this induces an action on the Grassmannian 
$\mathbb \GG(1, \PP(S^2V))$.
Given a non-degenerate pencil of quadrics in $\PP(V)$, the splitting type of the vector bundle $E$ determines a partition of the integer $r$ in $h$ parts, where $h=N-2r+1$. For every choice of constant rank $2r$ there are exactly $r$ possible sizes $N+1$ for these pencils, namely $N$ can vary from $2r$ to $3r-1$. On the other hand, if the rank and the order of the matrix are fixed, the number of parts $h$ of the partition of $r$ is determined. 

\smallskip

Our main result in this Section states that, for a fixed $r$, all possible values of $N$ are attained, and that the partitions of $r$ consisting of $h=N-2r+1$ parts completely characterize the orbits of pencils of quadrics of constant rank.

\smallskip

In our proof we will use the classification of the $\GL(N+1)$-orbits given in terms of minimal indices, see \cite[Chapter XII, \S 6]{Gantmacher}.

In fact, even if $\GL(N+1)$ acts on a pencil $\LL \subset \PP(S^2V)$ by congruence, one can also consider a different natural action of the general linear group on $\LL$, namely two pencils of matrices $a A+ bB$ and $\lambda L + \mu M$ are called {\em strictly equivalent} if there exist two non-singular matrices $P'$ and $P''$ such that $P'(a A + b B)P'' =\lambda L + \mu M$. The latter action implies the former if the matrices are symmetric or skew-symmetric \cite[Theorem 6, Chapter XII]{Gantmacher}: in particular, two pencils of quadrics are strictly equivalent  if and only if they are congruent.

\smallskip

Following the same notations as \cite{Gantmacher} (so slightly different than \cite{Fania_Mezzetti}), our construction is based on the following ``building blocks''.

\begin{defn}\label{def blocchi}
Let $r\ge 1$ be an integer, and $(r_1,\ldots,r_h)$ a partition of $r$, with $r_1  \le \ldots \le r_h$. Set $N = 2r + h -1$. Denote by $\LL_{(r_1,\ldots,r_h)}$ the pencil of $(N+1)\times (N+1)$ symmetric matrices of constant rank $2r$ constructed as follows. 

First, define the $r_i \times (r_i+1)$ matrix
\begin{equation}\label{M(ri)}
M_{r_i}:=\begin{pmatrix} 
a&b&0&0&\cdots&0\\
0&a&b&0&\cdots&0\\
\vdots&&\ddots&\ddots&&\vdots\\
0&\cdots&&0&a&b
\end{pmatrix},
\end{equation}
and the $(2r_i +1) \times (2r_i +1)$ symmetric block matrix
\begin{equation}\label{L(ri)}
\LL_{r_i}:=\left(\begin{array}{c;{2pt/2pt}c}
\renewcommand\arraystretch{2}
\:\: 0_{r,r}&M_{r_i} \:\:\:\:\:\\[-.8em]
&\\
\hdashline[2pt/2pt]
&\\[-.8em]
\:\:{}^t\!M_{r_i}&0_{r+1,r+1}
\end{array}\right).
\end{equation}

The pencil of quadrics $\LL_{(r_1,\ldots,r_h)}$ is the direct sum of the blocks $\LL_{r_i}$, so
\begin{equation}\label{L(r1,...,rh)}
\LL_{(r_1,\ldots,r_h)}:=
\left(\begin{array}{c|c|c|c}
\renewcommand\arraystretch{2}
\LL_{r_1}&&&\\[.4em]
\hline
&\LL_{r_2}&&\\[.4em]
\hline
&&\ddots &\\
\hline
&&&\LL_{r_h}
\end{array}\right),
\end{equation}
where all off-diagonal blank spaces are blocks of zeros.
\end{defn}

By combining the construction of the pencils $\LL_{(r_1,\ldots,r_h)}$ and the classification contained in  Theorem 7 and the subsequent remarks in \cite[Chapter XII, \S 6]{Gantmacher}, we obtain the following Theorem, that achieves a complete description of the $\GL(N+1)$-orbits of singular pencils of quadrics $\LL \subset \PP(S^2V)$ of constant rank.

\begin{thm}\label{partition}
Let $V$ be a complex vector space of dimension $N+1$, and let $\LL \subseteq \PP(S^2V)$ be a singular pencil of quadrics of constant rank $2r$. If $\LL$ is non-degenerate, it is $\GL(N + 1)$-equivalent by congruence and strict equivalence to a pencil of type $\LL_{(r_1,\ldots,r_h)}$ defined in \eqref{L(r1,...,rh)} for some partition $(r_1,\ldots,r_h)$ of $r$, with $r_1 \le  \ldots \le r_h$, $h=N+1-2r$, and whose associated vector bundle $E$ has splitting type precisely $(r_1, \ldots, r_h)$. 

Viceversa, for every integer $r \ge 1$ and every partition $(r_1, \ldots,r_h)$ of $r$, with $r_1 \le  \ldots \le r_h$, there exists a non-degenerate singular pencil of quadrics of constant rank $2r$ and size $N+1$, for all $2r \le N \le 3r-1$.
\end{thm}

\begin{rem} 
An alternative proof of Theorem \ref{partition} could be obtained by adapting to the symmetric case the proof of \cite[Theorem 3.12]{Fania_Mezzetti}, which is based on compression spaces and 1-generic matrices.
\end{rem}

\begin{rem}
If one wanted to take into consideration degenerate pencils, it would be enough to consider partitions of $r$ that admit 0 as a summand, with multiplicity corresponding to the number of copies of $\OP$ appearing in the splitting of the vector bundle $E$ in \eqref{succ esatta lunga}.
\end{rem}

To conclude this Section, we underline the fact that the content of Theorem \ref{partition} was already known, even though the relation between the classification of the orbits of singular pencils of quadrics of constant rank and the splitting type of the vector bundle has never been explicitly written down. In \cite{Dimca} the Author provides a geometric classification of the orbits, but the relation with the vector bundles is not clarified; on the other hand, in the recent work \cite{FJV} there is an explicit description of the splitting type of the bundles, but the Authors are interested in different properties than the orbits of pencils in the Grassmannian.

%%%%%%%%%%%%%%%%%%%%%%%%%%%%%%%%%%%%%%%%%%%%%%%%%%%%%%%%%%%%%%%%%%%%%%%%
%%%%%%%%%%%%%%%%%%%%%%%%%%%%%%%%%%%%%%%%%%%%%%%%%%%%%%%%%%%%%%%%%%%%%%%%
\section{Geometry of pencils of quadrics and their orbits}\label{examples}
%%%%%%%%%%%%%%%%%%%%%%%%%%%%%%%%%%%%%%%%%%%%%%%%%%%%%%%%%%%%%%%%%%%%%%%%
%%%%%%%%%%%%%%%%%%%%%%%%%%%%%%%%%%%%%%%%%%%%%%%%%%%%%%%%%%%%%%%%%%%%%%%%

We now want to study more in detail the geometry of %our 
pencils of quadrics of constant rank and their orbits. To this end, %we will sometime 
in this Section we use a different canonical form from the one given in Definition \ref{def blocchi} for the pencils with $h\ge 2$. It is analogous to the canonical form described in \cite{Fania_Mezzetti} in the skew-symmetric case, and is more convenient to understand the geometry of our pencils because it highlights that they are compression spaces. 

\smallskip

We start with some examples, describing the first cases where $r=1,2$ and $3$.

%%%%%%%%%%%%%%%%%%%%%%%%%%%%%%%%%%%%%%%%%%%%%%%%%%%%%%%%%%%%%%%%%%%%%%%%

\begin{ex}\label{r=1} 
The first (and easiest) example is $r=1$: then the only possible value for $N$ is $2$, and the only partition of $r$ is $(1)$, so there is a unique orbit, whose representative is the compression space
\begin{equation}\label{matrice O(1)}
\LL_{(1)}=
%\left(\begin{array}{c;{2pt/2pt}cc}
%0&a&b\\
%\hdashline[2pt/2pt]
%a&&\\
%b&&
%\end{array}\right)
\left(\begin{array}{c;{2pt/2pt}c}
0&{\begin{array}{cc}a&b\end{array}}\\
\hdashline[2pt/2pt]
{\begin{array}{c}
     a\\
     b 
\end{array}}&0_{2,2}
\end{array}\right).
\end{equation}
%This also agrees with results from \cite{Eisenbud_Harris}. 
The cokernel bundle from the exact sequence \eqref{succ esatta lunga} is $E=\OP(1)$. This is a pencil of conics in $\PP^2$, generated by $A=\{x_0x_1=0\}$ and $B=\{x_0x_2=0\}$, that split into a common line $S=\{x_0=0\}$ and a second line that goes through the point $P=[1:0:0]$. The base locus of the pencil is exactly the union of the line $S$, and the isolated point $P$. Notice that $S$ is swept by the singular points of the conics of the pencil.

The pencils belonging to the orbit of $\LL_{(1)}$ in the Grassmannian $\mathbb G(1,\PP^5)$ are determined by their base locus, that varies in the open subset of $\mathbb P^2\times {\mathbb{P}^2}^*$ of disjoint pairs point-line. Therefore the orbit in $\mathbb G(1,\PP^5)$ has dimension $4$.
\end{ex}

%%%%%%%%%%%%%%%%%%%%%%%%%%%%%%%%%%%%%%%%%%%%%%%%%%%%%%%%%%%%%%%%%%%%%%%%

\begin{ex}\label{r=2}
When $r=2$, the possible values of $N$ are $4$ and $5$, corresponding to the two partitions $(2)$ and $(1,1)$.

The first case gives a $5 \times 5$ symmetric matrix of constant rank 4:
\[\LL_{(2)}=
%\left(\begin{array}{cc;{2pt/2pt}ccc}
%&&a&b&0\\
%&&0&a&b\\
%\hdashline[2pt/2pt]
%a&0&&&\\
%b&a&&&\\
%0&b&&&
%\end{array}\right),
\left(\begin{array}{c;{2pt/2pt}c}
0_{2,2}&{\begin{array}{ccc}a&b&0\\0&a&b\end{array}}\\
\hdashline[2pt/2pt]
{\begin{array}{cc}
     a&0\\
     b&a\\
     0&b
\end{array}}&0_{3,3}
\end{array}\right),
\]
with associated line bundle $\OP(2)$.
The pencil is generated by $A=\{x_0x_2+x_1x_3\}$ and $B=\{x_0x_3+x_1x_4)\}$; its elements are cones over quadrics in $\PP^3$, having a single point as vertex. As the cones vary, their vertices describe a conic $\Gamma$ in the plane $S=\{x_0=x_1=0\}$.
The base locus is the union of the plane $S$ and the rational normal scroll of degree 3 in $\PP^4$ defined by the $2 \times 2$ minors of the matrix
\[\begin{pmatrix}
x_0&x_3&x_4\\
-x_1&x_2&x_3
\end{pmatrix}.\] 
The singular locus of the base locus is the conic $\Gamma$, which coincides with the improper intersection of the 2 irreducible components.

A pencil in this orbit is completely determined by its base locus, that is the union of  a rational normal scroll and a plane generated by a unisecant conic. From \cite{ellingsrud} we learn that the Hilbert scheme of these rational normal scrolls has dimension $12$; moreover the linear system of unisecant conics on such a surface has dimension $2$; it follows that the orbit has dimension $14$.
\smallskip

The partition $(1,1)$ of $r=2$ gives a $6\times 6$ symmetric matrix of constant rank 4, whose associated bundle is $E=\OP(1)\oplus \OP(1)$. As we mentioned at the beginning of the Section, we consider the following canonical form (here and in the next examples the blank spaces all represent zeros): 
%\[\LL_{(1,1)}=\left(\begin{array}{c|c}
%\LL_{1}& \\
%\hline
%&\LL_{1}
%\end{array}\right)=
%\left(\begin{array}{c|c}
%\begin{array}{c;{2pt/2pt}cc}
%&a&b\\
%\hdashline[2pt/2pt]
%a&&\\
%b&&
%\end{array}& \\
%\hline
%&\begin{array}{c;{2pt/2pt}cc}
%&a&b\\
%\hdashline[2pt/2pt]
%a&&\\
%b&&
%\end{array}
%\end{array}\right).
%\]
\[\tilde{\LL}_{(1,1)}=\left(\begin{array}{cc|cccc}
&&a&b&0&0\\
&&0&0&a&b\\
\hline
a&0&&&&\\
b&0&&&&\\
0&a&&&&\\
0&b&&&&
\end{array}\right).\]
Of course, $\tilde{\LL}_{(1,1)}$ is strictly equivalent to the block construction from Definition \ref{def blocchi}, namely:
\[\LL_{(1,1)}=\left(\begin{array}{c|c}
\LL_{1}& \\
\hline
&\LL_{1}
\end{array}\right).\]

Since the co-rank is $2$, the cones of this pencil have a line as vertex. The generators are $A=\{x_0x_2+x_1x_4\}$ and $B=\{x_0x_3+x_1x_5\}$, the base locus is reducible, and its components are the 3-dimensional linear space $S=\{x_0=x_1=0\}$ and a rational normal $3$-fold scroll of degree 3 in $\PP^5$, defined by the $2 \times 2$ minors of
\[\begin{pmatrix}
x_0&x_4&x_5\\
-x_1&x_2&x_3
\end{pmatrix}.\] 
The locus swept by vertices is a smooth quadric surface in $S$. By a  count of parameters similar to previous case, the dimension of the orbit is $26$: indeed, the dimension of the Hilbert scheme of rational normal cubic scrolls in $\mathbb P^5$ is $24$ and the linear system of unisecant quadrics  has dimension $2$.

One of the advantages of using the form $\tilde{\LL}_{(1,1)}$ lies precisely in the fact that the codimension 2 linear space $S$ contained in the base locus is now apparent, since we are dealing with a compression space. This phenomenon will generalize in the next cases.
\end{ex}

%%%%%%%%%%%%%%%%%%%%%%%%%%%%%%%%%%%%%%%%%%%%%%%%%%%%%%%%%%%%%%%%%%%%%%%%

\begin{ex}\label{r=3}
As a last series of examples, aiming to illustrate the general case,  we now consider the possible partitions of $r=3$. One has three possible values $6 \le N \le 8$, corresponding to the three partitions $(3)$, $(1,2)$ and $(1,1,1)$. By now we know that the representatives of their orbits are, respectively,
\[\LL_{(3)}=
\left(\begin{array}{c;{2pt/2pt}c}
0_{3,3}&{\begin{array}{cccc}a&b&0&0\\0&a&b&0\\0&0&a&b\end{array}}\\
\hdashline[2pt/2pt]
{\begin{array}{ccc}
     a&0&0\\
     b&a&0\\
     0&b&a\\
     0&0&b
\end{array}}&0_{4,4}
\end{array}\right), \quad \LL_{(1,2)}, \quad \hbox{and} \quad \LL_{(1,1,1)}.\]
%\[\LL_{(1,2)}=\left(\begin{array}{c|c}
%\LL_{1}& \\
%\hline
%&\LL_{2}
%\end{array}\right)=\left(\begin{array}{c|c}
%\begin{array}{c;{2pt/2pt}cc}
%&a&b\\
%\hdashline[2pt/2pt]
%a&&\\
%b&&
%\end{array}& \\
%\hline
%&{\begin{array}{cc;{2pt/2pt}ccc}
%&&a&b&0\\
%&&0&a&b\\
%\hdashline[2pt/2pt]
%a&0&&&\\
%b&a&&&\\
%0&b&&&
%\end{array}}
%\end{array}\right),\]
%\[\hbox{and} \quad \LL_{(1,1,1)}=\left(\begin{array}{c|c|c}
%\LL_{1}& &\\
%\hline
%&\LL_{1}&\\
%\hline
%&& \LL_{1}
%\end{array}\right)
%=\left(\begin{array}{c|c|c}
%\begin{array}{c;{2pt/2pt}cc}
%&a&b\\
%\hdashline[2pt/2pt]
%a&&\\
%b&&
%\end{array}& &\\
%\hline
%&\begin{array}{c;{2pt/2pt}cc}
%&a&b\\
%\hdashline[2pt/2pt]
%a&&\\
%b&&
%\end{array}&\\
%\hline
%&&\begin{array}{c;{2pt/2pt}cc}
%&a&b\\
%\hdashline[2pt/2pt]
%a&&\\
%b&&
%\end{array}
%\end{array}\right).\]
The base locus of the pencil $\LL_{(3)}$ in $\mathbb P^6$ is an irreducible quartic, complete intersection of the two quadrics $A=\{x_0x_3+x_1x_4+x_2x_5\}$ and $B=\{x_0x_4+x_1x_5+x_2x_6\}$; it is singular along a twisted cubic $C$ swept by the vertices and it contains the 3-dimensional linear space $S=\{ x_0=x_1=x_2=0\}$ spanned by $C$.  

\smallskip

To analyze the other two cases, we will again look at representatives that are strictly equivalent to $\LL_{(1,2)}$ and $\LL_{(1,1,1)}$, namely:
\[\tilde{\LL}_{(1,2)}=\left(\begin{array}{ccc|ccccc}
     &&&a&b&0&0&0\\ 
     &&&0&0&a&b&0\\ 
     &&&0&0&0&a&b\\ 
     \hline
     a&0&0&&&&&\\ 
     b&0&0&&&&&\\ 
     0&a&0&&&&&\\ 
     0&b&a&&&&&\\ 
     0&0&b&&&&& 
\end{array}
\right)\quad \hbox{and} \quad \tilde{\LL}_{(1,1,1)}=\left(\begin{array}{ccc|cccccc}
     &&&a&b&0&0&0&0\\ 
     &&&0&0&a&b&0&0\\ 
     &&&0&0&0&0&a&b\\ 
     \hline
     a&0&0&&&&&&\\ 
     b&0&0&&&&&&\\ 
     0&a&0&&&&&&\\ 
     0&b&0&&&&&&\\ 
     0&0&a&&&&&&\\
     0&0&b&&&&&&
\end{array}
\right).\]
Considering the kernels of these matrices, we easily see that in both cases  the Jacobian locus of the pencil is contained in the linear space $S=\{x_0=x_1=x_2=0\}$ of codimension 3 (so of dimension 4 and 5 respectively). The base locus  is irreducible in both cases and it is singular along the Jacobian locus, that is a rational normal scroll in $S$, $\PP(\OP(1) \oplus \OP(2))$ and $\PP(\OP(1) \oplus \OP(1) \oplus \OP(1))$ respectively. 
\end{ex}

We now describe the general case of a pencil $\LL=\LL_{(r_1,\ldots,r_h)}$ of constant rank $2r$ in $\mathbb P^N$, corresponding to the partition $(r_1,\cdots,r_h)$ of $r$, $h=N+1-2r$. Recall that we can write our $\LL$ as $\{aA + bB \ |\ [a:b] \in \PP^1\}$. We denote by $B(\LL)=A \cap B$ the base locus of $\LL$. It is a known fact that its singular locus is contained in the Jacobian locus $J(\LL)$ of $\LL$, the union of the vertices of the quadrics in the pencil, and such vertices are linear spaces of dimension $N-2r$.

\smallskip

As we did in the previous examples, we use a canonical form for the pencils that is slightly different from \eqref{L(r1,...,rh)}, and instead agrees with the notations used in \cite{Fania_Mezzetti}: given the $r_i \times (r_i+1)$ block $M_{r_i}$ defined in \eqref{M(ri)}, we set
\begin{equation}
\tilde{\LL}_{(r_1,\ldots,r_h)}:=
\left(\begin{array}{c|c}
     & {\begin{array}{c|c|c|c}
\renewcommand\arraystretch{2}
M_{r_1}&&&\\[.4em]
\hline
&M_{r_2}&&\\[.4em]
\hline
&&\ddots &\\
\hline
&&&M_{r_h}
\end{array}} \\
\hline
{\begin{array}{c|c|c|c}
\renewcommand\arraystretch{2}
{}^t\!M_{r_1}&&&\\[.4em]
\hline
&{}^t\!M_{r_2}&&\\[.4em]
\hline
&&\ddots &\\
\hline
&&&{}^t\!M_{r_h}
\end{array}}
     & 
\end{array}\right),
\end{equation}

where again the blank spaces have blocks of zeros.

From this canonical form, it is immediate to see that all these pencils correspond to compression spaces, because the associated matrices have a block of zeros of dimension $N+1-r$; a direct consequence is that the Jacobian locus $J(\LL)$ is contained in the linear space $S$ of dimension  $N-r$ defined by the equations $x_0=x_1=\cdots=x_{r-1}=0$.

Moreover, one easily computes that the Jacobian locus coincides with the singular locus of $B(\LL)$, which is irreducible, and it is exactly a rational normal scroll $\PP(\OP(r_1) \oplus \ldots \oplus \OP(r_h))$. 
Any element of the pencil is a cone over a smooth quadric of dimension $2r-2$, so it admits two families of linear spaces of dimension $(r-1)+(N-2r)+1=N-r$. Two quadrics of the pencil share a maximal linear subspace $S$ of dimension $N-r$ belonging to one of the two families, and are tangent along a rational normal scroll of type $r_1,\ldots, r_h$ in $S$.

\medskip

As a last remark ending this Section, we quote the article \cite{Segre}, a continuation and completion of the thesis of Corrado Segre, where he  studied the geometry of singular pencils of quadrics in $\PP^N$ of rank at most $k$, that he calls ``coni quadrici di specie $N-k$'', relating them to rational normal scrolls contained in their Jacobian locus.

%%%%%%%%%%%%%%%%%%%%%%%%%%%%%%%%%%%%%%%%%%%%%%%%%%%%%%%%%%%%%%%%%%%%%%%%
%%%%%%%%%%%%%%%%%%%%%%%%%%%%%%%%%%%%%%%%%%%%%%%%%%%%%%%%%%%%%%%%%%%%%%%%
\section{Orbits' dimension}\label{dimensione orbite}
%%%%%%%%%%%%%%%%%%%%%%%%%%%%%%%%%%%%%%%%%%%%%%%%%%%%%%%%%%%%%%%%%%%%%%%%
%%%%%%%%%%%%%%%%%%%%%%%%%%%%%%%%%%%%%%%%%%%%%%%%%%%%%%%%%%%%%%%%%%%%%%%%

 We recalled in Section \ref{first results} that the natural action of the group $\GL(N + 1)$ on $V=\C^{N+1}$ extends to the congruence action on $\PP(S^2V)$, and hence on the lines contained in $\PP(S^2V)$. Looking at pencils of quadrics as points in the Grassmannian $\GG(1,\PP(S^2V))$, we get  an action of $\GL(N+1)$ on the Grassmannian.
 We are interested in the orbits of singular pencils of quadrics $\LL \subseteq \PP(S^2V)$ of constant rank $2r$  under this latter action.
As we saw in Theorem \ref{partition}  non-degenerate pencils of quadrics  in $\mathbb P^N$ of constant rank $2r$ exist if and only if $2r\le N\le 3r-1$ and the orbits of these pencils correspond bijectively to the partitions $(r_1,\ldots,r_h)$ of $r$, with $1\le r_1 \le r_2 \le \ldots r_h$, where $h=N+1-2r$. 

This last Section contains our main result Theorem \ref{main}, namely we compute the dimension of all the orbits of pencils of singular quadrics of constant rank.
More precisely, for every partition $(r_1,\ldots,r_h)$ we describe explicitly the Lie algebra of the stabilizer of the pencil $\LL_{(r_1,\ldots,r_h)}$. 

\begin{thm}\label{main}
Let $r\ge 1$ be an integer, and $(r_1,\ldots,r_h)$ a partition of $r$, with $r_1  \le \ldots \le r_h$. Set $N = 2r + h -1$. Under the natural action of $\GL(N+1)$, the dimension of the stabilizer of the singular pencil $\LL_{(r_1,\ldots,r_h)}$ of symmetric matrices of size $N+1$ and constant rank $2r$ is
\begin{equation}\label{dimensione}
\delta(r_1,\ldots,r_h):=h + 4 + \sum_{i < j} (2r_j+1) + \#\{(i,j) \ | \ r_i=r_j\}.
\end{equation}
\end{thm}

\begin{cor}
The $\GL(N+1)$-orbits of singular pencils $\LL_{(r_1,\ldots,r_h)}$ of symmetric matrices of size $N+1$ and constant rank $2r$ have (affine) dimension $(N+1)^2-\delta(r_1,\ldots,r_h)$.
\end{cor}

The plan of the proof of Theorem \ref{main} is the following: we first analyze the case of partitions with only one part, i.e. pencils of symmetric matrices of constant corank $1$; we then consider the case of partitions with two parts, i.e. pencils of constant corank $2$. We obtain a complete description of the Lie algebra of the stabilizer in both cases. The key remark is then that, in the general case, due to the particular canonical form of the representatives of the orbits under consideration, a matrix $X$ in the Lie algebra of the stabilizer can be interpreted as a block matrix, where the blocks involved already appear and are described in the first two cases.

\smallskip

The next Lemma is probably well known. We report it here for completeness and because it is a fundamental ingredient for computing the Lie algebras of the stabilizers in the two cases $h=1,2$. 
\begin{lem}\label{lemma}
Let $\LL$ be the pencil  generated by the symmetric matrices $A$ and $B$, let $X$ be a $(N+1) \times (N+1)$ matrix with entries in $\C$. Then $X$ belongs to the Lie algebra of the stabilizer of $\LL$ for the action of $\GL(N+1)$ on the Grassmannian if and only if the following relations hold:
\begin{equation}\label{equazioni stab}
({}^t\!X A + A X) \wedge A \wedge B \: = \: ({}^t\!X B + B X) \wedge A \wedge B\:  = \: 0.
\end{equation}
\end{lem}
\begin{proof}
The point in the Grassmannian $\GG(1,\PP(S^2V))$ corresponding to the pencil $\LL$ via the Pl\" ucker map is $[A\wedge B]$. Its $\GL(N+1)$-orbit is the image of the map $\GL(N+1)\to \GG(1,\PP(S^2V))$ given by $X\mapsto ({}^t\!XAX)\wedge({}^t\!XBX)$. So the condition for  $X$ to belong to the stabilizer of $\LL$ is $[A\wedge B]=[({}^t\!XAX)\wedge({}^t\!XBX)]$. This is equivalent to the equations $({}^t\!XAX)\wedge A\wedge B=({}^t\!XBX)\wedge A\wedge B=0$.  Differentiating these equations at the origin we get the thesis.
\end{proof}

\begin{rem}
In the article \cite{DKS}, the Authors are interested in the same problem of computing the dimensions of orbits of pencils of symmetric matrices. But instead of interpreting them as points in the appropriate Grassmannian, they work with pairs of matrices generating the pencil, thus obtaining a different result from ours.
\end{rem}

We start with the partition having only $h=1$ part. We have a pencil of symmetric matrices of size $N+1=2r+1$ and rank $2r$, whose cokernel is the line bundle $E=\OP(r)$; the orbit representative is $\LL_{(r)}$, that we write in the following form, suitable to apply Lemma \ref{lemma}:

\begin{equation}\label{aA+bB per r}
    \LL_{(r)}=aA+bB=\left(\begin{array}{c;{2pt/2pt}c}
0_{r,r}&\begin{array}{ccccc} 
a&b&&&\\
&a&b&&\\
&&\ddots&\ddots&\\
&&&a&b
\end{array}\\
\hdashline[2pt/2pt]
\begin{array}{cccc} 
a&&&\\
b&a&\\
&b&\ddots&\\
&&\ddots&a\\
&&&b
\end{array}&0_{r+1,r+1}
\end{array}\right).
\end{equation}

\begin{prop}\label{dimensione O(r)}
Let $r \ge 1$ be an integer. The $\GL(2r+1)$-orbit of pencils of singular quadrics of constant rank $2r$ and order $2r+1$ has a stabilizer of dimension 5. 
The Lie algebra of the stabilizer is the vector space of matrices $X$ of the form:
\begin{equation}\label{matrice X}
X=\left(\begin{array}{c;{2pt/2pt}c}
\renewcommand\arraystretch{2}
\:\: X_1&0_{r,r+1} \:\:\:\:\:\\[-.8em]
&\\
\hdashline[2pt/2pt]
&\\[-.8em]
\:\:0_{r+1,r}&X_2
\end{array}\right),
\end{equation}
where:
\begin{enumerate}
    \item $X_1$ and $X_2$ are square matrices of order $r$ and $r+1$ respectively;
    \item both $X_1$ and $X_2$ are tridiagonal, i.e. all the elements out of the main diagonal, the sub-diagonal (the first diagonal below this), and the supradiagonal (the first diagonal above the main diagonal) are zero;
    \item the sub-diagonal, main diagonal, and supradiagonal of $X_1$ are respectively: $$y(r-1,r-2,\ldots,1),\  x_{00}(1,0,-1,-2,\ldots,-(r-2))+x_{11}(0,1,2,\ldots,r-1), \  z(1,2,\ldots,r-1);$$    
    \item the sub-diagonal, main diagonal, and supradiagonal of $X_2$ are respectively: $$-z(1,2,\ldots,r), \  (x_{00}-x_{11})(0,1,\ldots,r)+x_{33}(1,1,\ldots,1),\ -y(r,r-1,\ldots,1),$$
\end{enumerate} 
where $x_{00}, x_{11}, x_{33}, y,z $ are independent parameters.

\end{prop}
For instance, if $r=3$, $X$ is as follows:
\[\left(\begin{array}{ccc;{2pt/2pt}cccc}
x_{00}&z&0&&&&\\
2y&x_{11}&2z&&&&\\
0&y&2x_{11}-x_{00}&&&&\\
\hdashline[2pt/2pt]
&&&x_{33}&-3y&0&0\\
&&&-z&x_{00}-x_{11}+x_{33}&-2y&0\\
&&&0&-2z&2x_{00}-2x_{11}+x_{33}&-y\\
&&&0&0&-3z&3x_{00}-3x_{11}+x_{33}
\end{array}\right).\]

\begin{proof}
Let $X=(x_{ij})_{i,j=0,\ldots,N}$ be a matrix of unknowns.  If  $A, B$ are the matrices introduced in (\ref{aA+bB per r}), the elements of indices $i\le j$ in the symmetric matrices ${}^t\! XA+AX$ and ${}^t\! XB+BX$ are as described below:

\begin{equation}\label{equazioni caso 1A}
    ({}^t\!XA+AX)_{ij}=
    \left\{\begin{array}{lcl}
    x_{j+r,i}+x_{i+r,j}&\ \hbox{if}\ &0\le i\le j\le r-1\\
    x_{j-r,i}+x_{i+r,j}&&0\le i\le r-1,\, r\le j\le 2r-1\\
    x_{i+r,2r}&&0\le i\le r-1,\, j=2r\\
    x_{j-r,i}+x_{i-r,j}&&r\le i\le j\le 2r-1\\
    x_{i-r,2r}&&r\le i\le 2r-1,\, j=2r\\
    0 &&i=j=2r
    \end{array}\right.
\end{equation}
\begin{equation}\label{equazioni caso 1B}
    ({}^t\!XB+BX)_{ij}=
    \left\{\begin{array}{lcl}
    x_{j+r+1,i}+x_{i+r+1,j}&\ \hbox{if}\ &0\le i\le j\le r-1\\
    x_{i+r+1,r}&&0\le i\le r-1,\, j=r\\
    x_{j-r-1,i}+x_{i+r+1,j}&&0\le i\le r-1,\, r+1\le j\le 2r\\
    0&&i=r=j\\
    x_{j-r-1,i}+x_{i-r-1,j}&&r+1\le i\le j\le 2r\\
    x_{j-r-1,r}&&i=r,\, r+1\le j\le 2r
    \end{array}\right.
\end{equation}

%\begin{table}[h]
%    \begin{tabular}{c|c|c}
% $({}^t\!XA+AX)_{ij}$   &range of indices&\\
%    \hline
%&$0\le i\le j\le r-1$&$x_{j+r,i}+x_{i+r,j}$\\
%\hline
%&$0\le i\le r-1, r\le j\le 2r-1$&$x_{j-r,i}+x_{i+r,j}$\\
%\hline
%&$0\le i\le r-1, j=2r$&$x_{i+r,2r}$\\
%\hline
%&$r\le i\le j\le 2r-1$&$x_{j-r,i}+x_{i-r,j}$\\
%\hline
%&$r\le i\le 2r-1, j=2r $&$x_{i-r,2r}$\\
%\hline
%&$i=j=2r$&$0$\\
%\hline
%$({}^t\!XB+BX)_{ij}$   &range of indices&\\
%\hline
%&$0\le i\le j\le r-1$&$x_{j+r+1,i}+x_{i+r+1,j}$\\
%\hline
%&$0\le i\le r-1, j=r$&$x_{i+r+1,r}$\\
%\hline
%
%&$0\le i\le r-1, r+1\le j\le 2r$&$x_{j-r-1,i}+x_{i+r+1,j}$\\
%\hline
%&$i=r=j$ &$0$\\
%\hline
%&$r+1\le i\le j\le 2r$&$x_{j-r-1,i}+x_{i-r-1,j}$\\
%\hline
%&$i=r, r+1\le j\le 2r$&$x_{j-r-1,r}$\\
%\end{tabular}
%\caption{}
%\label{equazioni caso 1}
%\end{table}

In view of Lemma \ref{lemma}, $X$ belongs to the Lie algebra of the stabilizer of the orbit of $\LL_{(r)}$ if and only if it satisfies the equations (\ref{equazioni stab}), that are equivalent to a series of equations in the entries of each of the two matrices ${}^t\! XA+AX$ and ${}^t\! XB+BX$, and precisely:
\begin{enumerate}
    \item[(i)] vanishing of the elements with equal indices;
    \item[(ii)] vanishing of the elements with indices $0\le i<j\le r-1$, $r\le i<j\le 2r$, $(i,i+r+2), \ldots, (i,2r)$ for $i=0,\ldots, r-2$, and $(i,r), \ldots, (i,i+r-1)$ for $i=1,\ldots,r-1;$
    \item[(iii)] elements with indices $(0,r), (1,r+1), \ldots,(r-1,2r-1)$ must be two by two equal;
    \item[(iv)] elements with indices $(0,r+1), (1,r+2), \ldots,(r-1,2r)$ must be two by two equal.
\end{enumerate}
Now, using \eqref{equazioni caso 1A} and \eqref{equazioni caso 1B} together with (i) we get $x_{0,r}=x_{1,r+1}=\cdots=x_{r-1,2r-1}=x_{r,0}=\cdots=x_{2r-1,r-1}=0,$ and also 
$x_{0,r+1}=x_{1,r+2}=\cdots=x_{r-2,2r-1}=x_{r+1,0}=\cdots=x_{2r,r-1}=0;$ note that in all these cases the difference of the indices is either $r$ or $r+1$.

From the vanishings just obtained and  those in (ii) whose indices differ by $1$,  we get $x_{0,r-1}=x_{1,r}=\cdots=x_{r-1,2r-2}=x_{r-1,0}=\cdots=x_{2r-2,r-1}=0,$ 
$x_{0,r+2}=x_{1,r+3}=\cdots=x_{r-2,2r}=x_{r+2,0}=\cdots=x_{2r,r-2}=0,$ and also $x_{2r,r}=x_{r-1,2r}=0.$ 

We continue in this way, considering relations in (ii) whose indices differ by $2$ and so on, until we get all the claimed vanishings in matrix (\ref{matrice X}) and moreover the following $2r$ equations:
$$x_{0,1}+x_{r+1,r}=x_{1,2}+x_{r+2,r+1}=\cdots=x_{r-1,r}+x_{2r,2r-1}=0,$$
and the symmetric ones
$$x_{1,0}+x_{r,r+1}=x_{2,1}+x_{r+1,r+2}=\cdots=x_{r,r-1}+x_{2r-1,2r}=0.$$ 

The relations in (iii) and (iv) impose  $2r-2$ conditions on the elements of the main diagonal of $X$, and $2r-2$ conditions on the elements of the subdiagonal and supradiagonal of $X$, and precisely:
$$x_{0,0}+x_{r,r}=x_{1,1}+x_{r+1,r+1}=\cdots=x_{r-1,r-1}+x_{2r-1,2r-1},$$
$$x_{0,0}+x_{r+1,r+1}=x_{1,1}+x_{r+2,r+2}=\cdots=x_{r-1,r-1}+x_{2r,2r},$$
$$x_{1,0}+x_{r,r+1}=x_{2,1}+x_{r+1,r+2}=\cdots=x_{2r-1,2r},$$
$$x_{r+1,r}=x_{0,1}+x_{r+2,r+1}=\cdots=x_{r-2,r-1}+x_{2r,2r-1}.$$
Combining everything, we obtain for $X$ the expression in (\ref{matrice X}), with $z=x_{0,1}$ and $y=x_{r-1,r-2}$; the Proposition is proved.
\end{proof}

\medskip
Our description of the stabilizer compared with the known classification of Lie algebras of small dimension (\cite{KN}) gives the following result.

\begin{prop}\label{descrizione algebra}
The Lie algebra of the stabilizer of the $\GL(2r+1)$-orbit of pencils of quadrics of constant rank $2r$ and order $2r+1$ described in Proposition \ref{dimensione O(r)} is isomorphic to $\mathfrak{sl}_2 \ltimes \C^2$.
\end{prop}

\begin{proof}
From the detailed description of the Lie algebra of the stabilizer given in Proposition \ref{dimensione O(r)}, one sees that its elements depend on $5$ independent parameters, namely any element $X$ in this Lie algebra is $X=X(x_{00},x_{11},x_{33},y,z)$. With obvious notation, let us call
\[\mathcal{C}_1=X(1,0,0,0,0), \quad \mathcal{C}_2=X(0,1,0,0,0),\]
\[\mathcal{X}=X(0,0,0,0,1), \quad \mathcal{Y}=X(0,0,0,1,0), \quad \mathcal{Z}=X(r-1, r-3, -r, 0,0).\]

If we compute the bracket of these elements, we get that $[\mathcal{C}_1,\mathcal{C}_2]=0$ and 
\[\begin{cases}
[\mathcal{X},\mathcal{Y}]=\mathcal{Z}\\
[\mathcal{Z},\mathcal{X}]=2\mathcal{X}\\
[\mathcal{Z},\mathcal{Y}]=-2\mathcal{Y}
\end{cases}\]
which tells us that $\C^2=<\mathcal{C}_1,\mathcal{C}_2>$ and $\mathfrak{sl}_2=<\mathcal{X},\mathcal{Y},\mathcal{Z}>$. The fact that 
\[\begin{cases}
[\mathcal{C}_1,\mathcal{X}]=\mathcal{X}=-[\mathcal{C}_2,\mathcal{X}]\\
[\mathcal{C}_1,\mathcal{Y}]=-\mathcal{Y}=-[\mathcal{C}_2,\mathcal{Y}]\\
[\mathcal{C}_1,\mathcal{Z}]=[\mathcal{C}_2,\mathcal{Z}]=0
\end{cases}\]
allows us to conclude that our Lie algebra falls into the first case in the classification table appearing in \cite[Section 4]{KN}, namely the semidirect product $\mathfrak{sl}_2 \ltimes \C^2$.
\end{proof}

\medskip

When the partition has $h=2$ parts, the balanced and unbalanced case have two different behaviors, as explained in the following result.

\begin{prop}\label{dimensione O(k)+O(r-k)}
Let $r \ge 1$ be an integer. The $\GL(2r+2)$-orbit of pencils of singular quadrics of constant rank $2r$ and order $2r+2$, whose associated bundle is $\OP(r_1)\oplus \OP(r_2)$, with $r_1+r_2=r$, and $r_1 \le r_2$, has stabilizer of dimension 
\begin{enumerate}
    \item $2r_2+8=r+8$ when $r$ is even and $r_1=r_2=\frac{r}{2}$;
    \item $2r_2+7$ when $r_1 < r_2$.
\end{enumerate} 
\end{prop}

\begin{proof}
In the notation of Section \ref{first results}, a representative of the orbit is the matrix 
\[\LL_{(r_1,r_2)}=aA+bB= \left(\begin{array}{c|c}
\LL_{r_1}& \\
\hline
&\LL_{r_2}
\end{array}\right).\] 
We also introduce the notation $A=\left(\begin{array}{c|c}
A_1& \\
\hline
&A_{2}
\end{array}\right)$, $B=\left(\begin{array}{c|c}
B_1& \\
\hline
&B_{2}
\end{array}\right)$, where  $A_i, B_i$ are matrices of order $2r_i+1$, for $i=1,2$. 

Let $X=(x_{ij})_{i,j=0,\ldots,N}$ be a matrix of unknowns. We write $X$ as a block matrix as follows:
\[X=\left(\begin{array}{c|c}
\renewcommand\arraystretch{2}
X_{11}&X_{12}\\[-1.1em]
&\\
\hline
&\\[-1em]
X_{21}&X_{22}
\end{array}\right)=\left(\begin{array}{c|c}
\renewcommand\arraystretch{2}
(x_{ij})_{\substack{i=0,\ldots,2r_1+1 \\ j=0,\ldots,2r_1+1}}&(x_{ij})_{\substack{i=0,\ldots,2r_1+1 \\ j=2r_1+2, \ldots, N}}\\[-1.1em]
&\\
\hline
&\\[-1em]
(x_{ij})_{\substack{i=2r_1+2,\ldots,N \\ j=0,\ldots,2r_1+1}}&(x_{ij})_{\substack{i=2r_1+2,\ldots,N \\ j=2r_1+2,\ldots,N}}
\end{array}\right)\]
where $X_{ii}$ are square matrices of order $(2r_i+1)$,  and $X_{12}$, $X_{21}$ have order   $(2r_1+1)\times (2r_2+1)$ and $(2r_2+1)\times(2r_1+1)$ respectively. 

Then ${}^t\!XA+AX$ and ${}^t\!XB+BX$ can be written as block matrices as well, and precisely:
\begin{equation}
{}^t\!XA+AX=
\left(\begin{array}{c|c}
\renewcommand\arraystretch{2}
{}^t\!X_{11}A_1+A_1X_{11}&{}^t\!X_{21}A_2+A_1X_{12}\\[-1.1em]
&\\
\hline
&\\[-1em]
{}^t\!X_{12}A_1+A_2X_{21}&{}^t\!X_{22}A_2+A_2X_{22}
\end{array}\right),
\end{equation}
and similarly for $B$.
Lemma \ref{equazioni stab} implies that $X$ belongs to the Lie algebra of the stabilizer if and only if equations (\ref{equazioni stab}) are satisfied. We analyze separately what this means for the diagonal blocks $X_{11}, X_{22}$ and for the off-diagonal blocks $X_{12}, X_{21}$ of $X$.

\smallskip

\underline{Diagonal blocks.} We use Proposition \ref{dimensione O(r)}: $X_{11}$, $X_{22}$ must belong to the Lie algebras of the stabilizers of the orbits of $\LL_{(r_1)}$ and $\LL_{(r_2)}$ respectively, therefore each of them depends on $5$ parameters and has the form described in Proposition \ref{dimensione O(r)}. But equations (\ref{equazioni stab}) imply that the parameters appearing in $X_{11}$ and $X_{22}$ are not independent, and precisely, after fixing the $5$ parameters required to describe $X_{11}$, an explicit computation shows that only one new parameter is needed to describe $X_{22}$, therefore the two diagonal blocks depend on a total of $6$ parameters.

\smallskip

\underline{Off-diagonal blocks.} The matrices ${}^t\!X_{21}A_2+A_1X_{12}$ and ${}^t\!X_{12}A_1+A_2X_{21}$ are the transpose of each other, and they both have to be the zero matrix. The same holds for ${}^t\!X_{21}B_2+B_1X_{12}$ and ${}^t\!X_{12}B_1+B_2X_{21}$. 

From the explicit expressions of their entries, we get
 the following conditions:
\begin{equation}\label{prime condizioni}
  x_{a,b}+x_{i,j}=0 \ \text{for any}\ 2r_1+1\le a,j\le 2r,\ 0\le i,b\le 2r_1-1 \ \text{with} \   |b-i|=r_1, |a-j|=r_2,
\end{equation}
\begin{equation}\label{seconde condizioni}
  x_{a,b}+x_{i,j}=0 \ \text{for any}\ 2r_1+1\le a,j\le 2r+1,\ 0\le i,b\le 2r_1 \ \text{with} \   |b-i|=r_1+1, |a-j|=r_2+1.  
\end{equation}
We also get a first series of four vanishings, referring to the last and the central columns of $X_{12}$ and $X_{21}$:
\begin{enumerate}
    \item[(i)] the last column of $X_{12}$ except its last element:
    \begin{equation*}
    x_{0,2r+1}=x_{1,2r+1}=\cdots=x_{2r_1-1, 2r+1}=0,
      \end{equation*}
    \item[(ii)] the central column of $X_{12}$, of index $2r_1+r_2+1$, except its central element $x_{r_1, 2r_1+r_2+1}$;
     \item[(iii)] the last column of $X_{21}$ except its last element: 
     \begin{equation*}
    x_{2r_1+1,2r_1}=x_{2r_1+2,2r_1}=\cdots=x_{2r,2r_1}=0,
    \end{equation*}
    \item[(iv)] the central column of $X_{21}$, of index $r_1$, except its central element $x_{2r_1+r_2+1, r_1}.$
\end{enumerate}
The vanishing of these columns, together with conditions (\ref{prime condizioni}) and (\ref{seconde condizioni}), implies, in order, the following second series of vanishings, referring to the rows of the two matrices:
\begin{enumerate}
    \item[(i)] the row of index $2r_1+r_2$ of $x_{21}$, except the element
    $X_{2r_1+r_2,r_1-1};$ this is the row above the middle;
    \item[(ii)] the first row of $X_{21}$ except its first element $x_{2r_1+1, 0};$
    \item[(iii)] the row of index $r_1-1$ of $X_{12}$ except $x_{r_1-1,2r_1+r_2};$ this is the row above the middle;
    \item[(iv)] the first row of $X_{12}$ except its first element $x_{0, 2r_1+1}.$
\end{enumerate}

\smallskip

We now analyze separately the two cases {\em (1)} and {\em (2)} in our statement.

\smallskip

Case {\em(1)}: when $r_1=r_2$, $X_{12}, X_{21}$ are square matrices. Going on with the argument above, we deduce that in both $X_{12}$ and $X_{21}$ all the elements above the central row and to the right of the central column are zero, except those of the main diagonal.
Moreover, the first $r_2$ entries of the main diagonal of $X_{12}$ are equal to each other and also to the last $r_2$ elements of the main diagonal of $X_{21}$, and similarly the last $r_2$ elements of the main diagonal of $X_{12}$ are equal to each other and also to the first $r_2$ elements of the main diagonal of $X_{21}$. 

We are left to analyze the two rectangles of order $(r_2+1) \times r_2$ in the lower left corner: from conditions (\ref{prime condizioni}) and (\ref{seconde condizioni}) we get that they depend on $2r_2$ parameters, independent of those previously considered. More precisely, we can divide each of the two rectangles into its $2r_2$ anti-diagonals; each of them results to be formed by elements all equal to each other and to those of the same anti-diagonal of the other matrix. 

All in all, there are $2+2r_2$ independent parameters for this case {\em (1)}. For the reader's convenience, we illustrated the case $(2,2)$ in Figure \ref{figura O(2)+O(2)}.

\begin{figure}
\begin{center}
\includegraphics[width=14truecm]{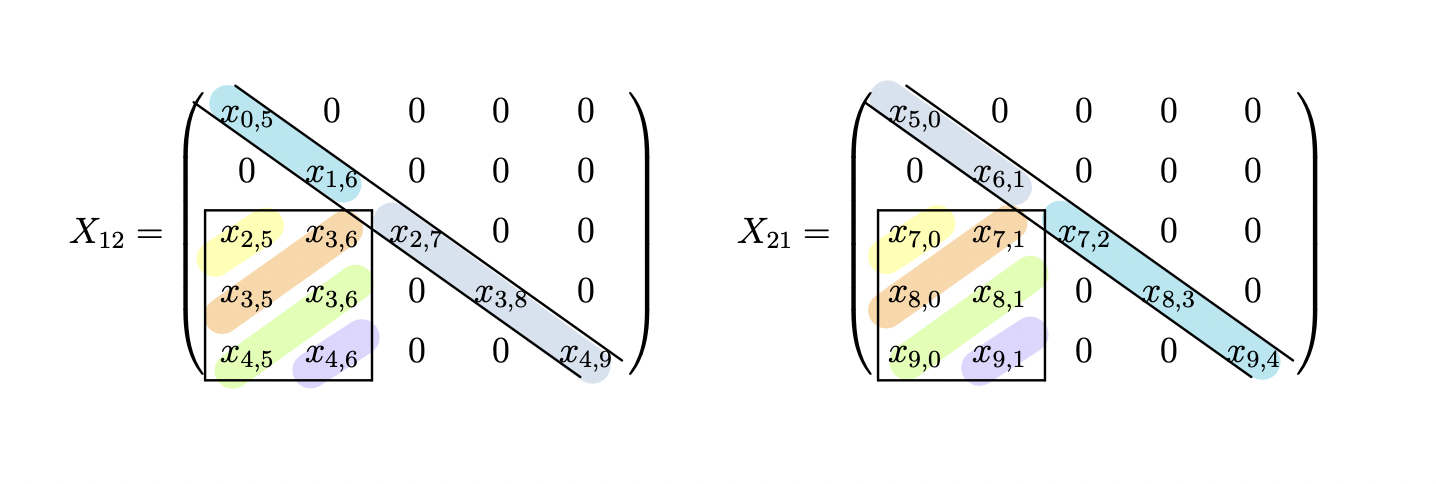}
\caption{Structure of the submatrices $X_{12}$ and $X_{21}$ in an element of the Lie algebra of the stabilizer of $\LL_{(2,2)}$: entries that are equal (up to a sign) are highlighted with the same color.}\label{figura O(2)+O(2)}
\end{center}
\end{figure}

\smallskip

Case {\em (2)}: assume now $r_1<r_2$. We obtain the vanishing of the entire first $r_1$ rows of $X_{12}$ and of the last $r_1+1$ columns of $X_{21}$. Now we need to look at the last $r_1+1$ rows of $X_{12}$ and the first $r_1$ columns of $X_{21}$. The former is divided into two blocks $\alpha_{12}$ and $\beta_{12}$ of size $(r_1+1) \times r_2$ and $(r_1+1) \times (r_2+1)$ respectively, while the latter is divided into two blocks $\alpha_{21}$ and $\beta_{21}$ of size $r_2 \times r_1$ and $(r_2+1)\times r_1$ respectively. All entries in each of the $r_2+2$ anti-diagonals of $\alpha_{12}$ are equal to each other, and the same is true for the $r_2+2$ anti-diagonals of $\beta_{21}$. Moreover, these diagonals are paired, in the sense that they depend in order exactly on the same $r_2+2$ parameters. Finally, the same relations hold for the $r_2-1$ principal diagonals of the blocks $\beta_{12}$ and $\alpha_{21}$, with the difference that this time all entries above and below these $r_2-1$ principal diagonals are zero. (By ``principal diagonal'' we mean a maximal length diagonal with $r_1+1$ entries in $\beta_{12}$ and $r_1$ entries in $\alpha_{21}$.)

All in all, there are $(r_2+2)+(r_2-1)=2r_2+1$ independent parameters for this case {\em (2)}. Figure \ref{figura O(2)+O(3)} illustrates the case $(2,3)$.

\begin{figure}
\begin{center}
\includegraphics[width=14truecm]{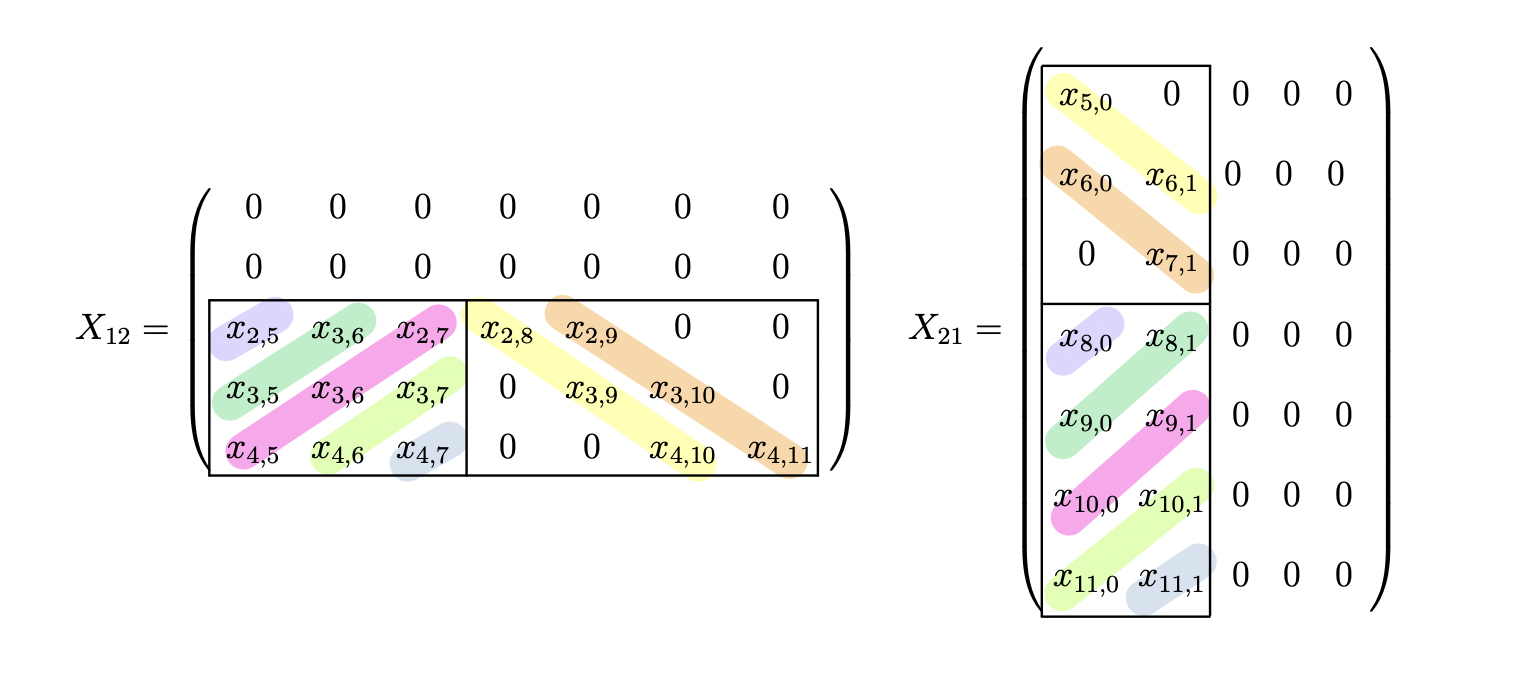}
\caption{Structure of the submatrices $X_{12}$ and $X_{21}$ in an element of the Lie algebra of the stabilizer of $\LL_{(2,3)}$: again, the entries that are equal (up to a sign) are highlighted with the same color.}\label{figura O(2)+O(3)}
\end{center}
\end{figure}

\smallskip

Notice that the unknowns appearing in the on and off-diagonal blocks are independent from each other: this means that we only need to add the number of independent parameters coming from the off-diagonal blocks to the $6$ ones needed for the diagonal blocks. This concludes the proof in both cases.
\end{proof}

\begin{proof}[Proof of Theorem \ref{main}]
We mimic and generalize the proof of Proposition \ref{dimensione O(k)+O(r-k)}. Given a pencil $\LL_{(r_1,\ldots,r_h)}$ in the canonical form \eqref{L(r1,...,rh)} and generated by $A$ and $B$, with obvious notation we write
\[A=
\left(\begin{array}{c|c|c|c}
\renewcommand\arraystretch{2}
A_{1}&&&\\[.4em]
\hline
&A_2&&\\[.4em]
\hline
&&\ddots &\\
\hline
&&&A_h
\end{array}\right)
\quad \hbox{and} \quad 
B=
\left(\begin{array}{c|c|c|c}
\renewcommand\arraystretch{2}
B_{1}&&&\\[.4em]
\hline
&B_2&&\\[.4em]
\hline
&&\ddots &\\
\hline
&&&B_h
\end{array}\right).\]

To describe the matrices $X$ belonging to the Lie algebra of the stabilizer of $\LL_{(r_1,\ldots,r_h)}$ we use Lemma \ref{lemma}. We write a general matrix of unknowns $X=(x_{ij})_{i,j=0,\ldots,N}$ as a block matrix with the same type of blocks $X_{ij}$ as above, each of size $(2r_i+1) \times (2r_j+1)$:
\[X=
\left(\begin{array}{c|c|c|c}
\renewcommand\arraystretch{2}
X_{11}&X_{12}&\ldots&X_{1h}\\[.4em]
\hline
X_{21}&X_{22}&&\vdots\\[.4em]
\hline
\vdots&&\ddots &\vdots\\
\hline
X_{1h}&\ldots&\ldots&X_{hh}
\end{array}\right).
\]
Then ${}^t\!XA+AX$ can also be written as a block matrix, where the square blocks on the diagonal have the form
\[
{}^t\!X_{ii}A_i+A_iX_{ii},
\]
while the off-diagonal ones with $i< j$ are 
\[
{}^t\!X_{ji}A_j+A_iX_{ij},
\]
and similarly for $B$.
As in the proof of Proposition \ref{dimensione O(k)+O(r-k)}, the upper left diagonal block $X_{11}$ depends on 5 independent parameters, and each other diagonal block contributes with 1 more degree of freedom. This accounts for $5+(h-1)=4+h$ parameters.
The off-diagonal blocks $X_{ij}$ and its symmetric $X_{ji}$ are in the same relation described for $X_{12}$ and $X_{21}$ in the proof of Proposition \ref{dimensione O(k)+O(r-k)}, so each pair accounts for $2r_j+2$ if $r_i=r_j$, and $2r_j+1$ if $r_i < r_j$. 

Since the blocks $X_{ij}$ and $X_{k\ell}$ are independent for $(i,j) \neq (k, \ell)$, the total number of parameters is 
\[4+h + \sum_{i < j} (2r_j+1) + \#\{(i,j) \ | \ r_i=r_j\},\]
and this concludes our proof.
\end{proof}

To illustrate our result, we collected in Table \ref{tabella} all orbits of pencils of quadrics of constant rank $2r$, $r \le 6$, their dimension, and the dimension of their stabilizer.

\begin{table}
    \begin{tabular}{c|c|c|c|c|c}
    $r$&$h$&partition&$N=2r+h-1$&dim orbit&dim stabilizer\\
    \hline
1&1&(1)&2&4&5\\
\hline
2&1&(2)&4&20&5\\
&2&(1,1)&5&26&10\\
\hline
3&1&(3)&6&44&5\\
&2&(2,1)&7&53&11\\
&3&(1,1,1)&8&62&19\\
\hline
4&1&(4)&8&76&5\\
&2&(2,2)&9&88&12\\
&2&(1,3)&9&87&13\\
&3&(1,1,2)&10&100&21\\
&4&(1,1,1,1)&11&112&32\\
\hline
5&1&(5)&10&116&5\\
&2&(2,3)&11&131&13\\
&2&(1,4)&11&129&15\\
&3&(1,2,2)&12&146&23\\
&3&(1,1,3)&12&144&25\\
&4&(1,1,1,2)&13&161&35\\
&5&(1,1,1,1,1)&14&176&49\\
\hline
6&1&(5)&12&164&5\\
&2&(3,3)&13&182&14\\
&2&(2,4)&13&181&15\\
&2&(1,5)&13&179&16\\
&3&(2,2,2)&14&200&25\\
&3&(1,2,3)&14&199&26\\
&3&(1,1,4)&14&196&29\\
&4&(1,1,2,2)&15&218&38\\
&4&(1,1,1,3)&15&215&40\\
&5&(1,1,1,1,2)&16&236&53\\
&6&(1,1,1,1,1,1)&17&254&70\\
\end{tabular}
\caption{Dimension of orbits of pencils of quadrics and their stabilizers.}
\label{tabella}
\end{table}

\smallskip

Looking at Table \ref{tabella}, it is interesting to observe the phenomenon occurring when there are two different partitions of $r$ of the same length. As expected from the behaviour of a rational normal scroll $\PP(\OP(r_1) \oplus \OP(r_2))$ degenerating to a $\PP(\OP(r_1-1) \oplus \OP(r_2+1))$, the dimension of the relative orbit increases.

\thispagestyle{empty}
\bibliographystyle{amsalpha}
\bibliography{biblio_quadriche}

\providecommand{\bysame}{\leavevmode\hbox to3em{\hrulefill}\thinspace}
\providecommand{\MR}{\relax\ifhmode\unskip\space\fi MR }
\providecommand{\MRhref}[2]{%
  \href{http://www.ams.org/mathscinet-getitem?mr=#1}{#2}
}
\providecommand{\href}[2]{#2}
\begin{thebibliography}{BFM15}

\bibitem[BFM15]{bor-fan-mez-quadriche}
A.~Boralevi, M.L. Fania, and E.~Mezzetti, \emph{Quadric surfaces in the
  {P}faffian hypersurface in $\mathbb{P}^{14}$}, Linear and Multilinear Algebra
  (2015), Published online 05 Mar 2021.

\bibitem[BM15]{bo_me_piani}
A.~Boralevi and E.~Mezzetti, \emph{Planes of matrices of constant rank and
  globally generated vector bundles}, Ann. Inst. Fourier \textbf{65} (2015),
  no.~5, 2069--2089.

\bibitem[Dim83]{Dimca}
A.~Dimca, \emph{A geometric approach to the classification of pencils of
  quadrics}, Geom. Dedicata \textbf{14} (1983), no.~2, 105--111.

\bibitem[DKS14]{DKS}
A.~Dmytryshyn, B.~K{\aa}gstr\"{o}m, and V.V. Sergeichuk, \emph{Symmetric matrix
  pencils: codimension counts and the solution of a pair of matrix equations},
  Electron. J. Linear Algebra \textbf{27} (2014), 1--18.

\bibitem[Ell75]{ellingsrud}
G.~Ellingsrud, \emph{Sur le sch\' ema de {H}ilbert des vari\' et\' es de
  codimension $2$ dans $\mathbf{P}^e$ \`a cône de {C}ohen-{M}acaulay}, Ann.
  Sci. \' {E}cole Norm. Sup. (4) \textbf{8} (1975), no.~4, 423--431.

\bibitem[FJV21]{FJV}
D.~Faenzi, M.~Jardim, and J.~Vall\`es, \emph{Logarithmic sheaves of complete
  intersections}, arXiv:2106.14453, 2021.

\bibitem[FM11]{Fania_Mezzetti}
M.L. Fania and E.~Mezzetti, \emph{Vector spaces of skew-symmetric matrices of
  constant rank}, Linear Algebra Appl. \textbf{434} (2011), 2383--2403.

\bibitem[FMS21]{Fevola-Sturmfels}
C.~Fevola, Y.~Mandelshtam, and B.~Sturmfels, \emph{Pencils of quadrics: old and
  new}, Matematiche (Catania) \textbf{76} (2021), no.~2, 319--335.

\bibitem[Gan59]{Gantmacher}
F.~R. Gantmacher, \emph{The theory of matrices. {V}ols. 1, 2}, Chelsea
  Publishing Co., New York, 1959, Translated by K. A. Hirsch.

\bibitem[IL99]{Ilic_JM}
B.~Ilic and J.M. Landsberg, \emph{On symmetric degeneracy loci, spaces of
  symmetric matrices of constant rank and dual varieties}, Math. Ann.
  \textbf{314} (1999), no.~1, 159--174.

\bibitem[KN87]{KN}
A.A. Kirillov and Y.A. Neretin, \emph{The {V}ariety ${A}_n$ of
  $n$-{D}imensional {L}ie {A}lgebra {S}tructures}, Amer. Math. Soc Transl. (2)
  \textbf{137} (1987), 21--30.

\bibitem[MM05]{Manivel_Mezzetti}
L.~Manivel and E.~Mezzetti, \emph{On linear spaces of skew-symmetric matrices
  of constant rank}, Manuscripta Math. \textbf{117} (2005), no.~3, 319--331.

\bibitem[Seg84]{Segre}
C.~Segre, \emph{Ricerche sui fasci di coni quadrici in uno spazio lineare
  qualunque}, Atti R. Acc. Scienze Torino \textbf{19} (1883-84), 878--896,
  anche in {C}orrado {S}egre, {O}pere, a cura della {U}nione {M}atematica
  {I}taliana, {V}ol. {III}, {E}dizione {C}remonese, {R}oma, 1961, p. 485–501.

\bibitem[Tho91]{Thompson}
R.C. Thompson, \emph{Pencils of complex and real symmetric and skew matrices},
  Linear Algebra Appl. \textbf{147} (1991), 323--371.

\end{thebibliography}

\end{document}